\documentclass[12pt,a4paper]{article}

\usepackage[latin2]{inputenc}
\usepackage{amssymb}
\usepackage{paralist}
\usepackage{amsmath, calc}
\usepackage{amsfonts}
\usepackage{amsthm}
\usepackage{graphicx}
\usepackage{fullpage}
\usepackage{enumerate}
\usepackage{mathdesign}
\newtheorem{thm}{Theorem}

\newtheorem{prob}{Problem}

\newtheorem{lem}{Lemma}

\begin{document}

\title{Dense sumsets of Sidon sequences}
\author{S\'andor Z. Kiss \thanks{Institute of Mathematics, Budapest
University of Technology and Economics, H-1529 B.O. Box, Hungary;
ksandor@math.bme.hu;
This author was supported by the National Research, Development and Innovation Office NKFIH Grant No. K115288 and K129335. 
This paper was supported by the J\'anos Bolyai Research Scholarship of the Hungarian Academy of Sciences. Supported by the \'UNKP-20-5 New National Excellence Program 
of the Ministry for Innovation and Technology from the source of the National Research Development and Innovation Fund.}, Csaba
S\'andor \thanks{Institute of Mathematics, Budapest University of
Technology and Economics, MTA-BME Lend\"ulet Arithmetic Combinatorics Research Group H-1529 B.O. Box, Hungary, csandor@math.bme.hu.
This author was supported by the NKFIH Grants No. K129335. Research supported 
by the Lend\"ulet program of the Hungarian Academy of Sciences (MTA), under grant number LP2019-15/2019.} 
}
\date{}
\maketitle

\begin{abstract}
\noindent Let $k \ge 2$ be an integer. We say a set $A$ of positive integers is an asymptotic basis of order $k$ if every large enough positive integer can be represented as the sum of $k$ terms from $A$. A set of positive integers $A$ is
called Sidon set if all the two terms sums formed by the elements of $A$ 
are different. Many years ago P. Erd\H{o}s, A. S\'ark\"ozy and V. T. S\'os
asked whether there exists a Sidon set which is asymptotic basis of order $3$. 
In this paper we prove the existence of a Sidon set $A$ with positive lower 
density of the three fold sumset $A + A + A$ by using probabilistic methods. 


\end{abstract}
\section{Introduction}

Let $\mathbb{N}$ denote the set of positive integers. Let $h, k \ge 2$ be
integers. Let $A \subset \mathbb{N}$ be an infinite set and let 
$R_{h,A}(n)$ denote the number of solutions of the equation 
\[
a_{1} + a_{2} + \dots + a_{h} = n, \hspace*{3mm} a_{1} \in
A, \dots, a_{h} \in A, \hspace*{3mm} a_{1} \le
a_{2} \le \dots{} \le a_{h},
\]
\noindent where $n \in \mathbb{N}$. Put 
\[
A(n) = \sum_{\overset{a \in \mathcal{A}}{a \le n}}1.
\]
A set of positive integers $A$ is
called a $B_h[g]$ set if for every $n \in \mathbb{N}$, the number of 
representations of $n$ as the sum of $h$ terms in the above form is at most $g$,
that is $R_{h,A}(n) \le g$. Obviously, the Sidon sets are the $B_{2}[1]$ sets.   
We say a set $A \subset \mathbb{N}$ is an asymptotic basis
of order $k$ if there exists a positive integer $n_{0}$ such that
$R_{k,A}(n) > 0$ for $n > n_{0}$. Let $A + A + A = \{a + a^{'} + a^{''}: a,a^{'},a^{''}\in A\}$. Moreover, for sets $A$, $B$ of positive integers, we define the representation function
\[
R_{A+B}(n) = \{a + b: a\in A, b\in B\}.
\]
In \cite{ET} and \cite{ES}, P. Erd\H{o}s, A. S\'ark\"ozy and V. T. S\'os posed the following question.

\begin{prob}
Does there exist a Sidon set which is an asymptotic basis of order $3$?
\end{prob}
This problem is still unsolved and it seems to be very difficult. A few years ago, Cilleruelo proved the following weaker result \cite{CE}.

\begin{thm}[Cilleruelo, 2015]
There exists a $B_2[2]$ set which is an asymptotic basis of order $3$.
\end{thm}
Cilleruelo applied the probabilistic method to prove the existence of such a 
set.
It is easy to see that if $A$ is an asymptotic basis of order $3$, then $A(x) \ge (\sqrt[3]{6} + o(1))\sqrt[3]{x}$, where $A(x)$ denotes the number of elements of the set of positive integers $A$ up to $x$.
By using the greedy algorithm \cite{CL}, one can construct a Sidon set $A$ with $A(x) \gtrapprox  (\sqrt[3]{2} + o(1))\sqrt[3]{x}$. 

It seems to be very difficult to construct denser infinite Sidon sets. In 1981, Ajtai, Koml\'os and Szemer\'edi \cite{AJK} proved the existence of an infinite Sidon set with $A(n) \gg \sqrt{n\log n}$ for every large enough $n$. Later, Ruzsa
\cite{RI} improved on this result to $A(n) \gg n^{\sqrt{2}-1+o(1)}$, whenever 
$n \rightarrow \infty$, by using probabilistic arguments.
Constructing directly an asymptotic basis of order $3$, which is a Sidon set 
seems to be hopeless as well. The above mentioned constructions do not help 
to settle Problem 1.

Furthermore, it is easy to see that a Sidon set cannot be an asymptotic basis of
order $2$. J. M. Deshouillers and A. Plagne
\cite{DE} constructed a Sidon set which is an asymptotic basis of order at most
$7$. In \cite{SB}, it was proved the existence of Sidon sets which are asymptotic bases of order $5$ by using probabilistic methods. 
It was proved by Cilleruelo \cite{CE} and independently at the same time 
by Kiss, Rozgonyi and S\'andor \cite{KR} that

\begin{thm}[Cilleruelo, Kiss-Rozgonyi-S\'andor, 2013]
There exists a Sidon set which is an asymptotic basis of order $4$.
\end{thm}
Moreover, Cilleruelo \cite{CE} proved the following, a little bit stronger 
result. A set $A$ is called an asymptotic basis of order $3 + \varepsilon$ if every sufficiently large integer $n$ can be written as a sum of $4$ terms from $A$ such that one of the terms is less than $n^{\varepsilon}$. 

\begin{thm}[Cilleruelo, 2015]
For every $\varepsilon > 0$, there exists a Sidon set which is an asymptotic 
basis of order $3 + \varepsilon$.
\end{thm}

Define the lower asymptotic density of a set of natural numbers by
\[
\underline{d}(A) = \liminf_{n \rightarrow \infty}\frac{A(n)}{n},
\]
and the upper asymptotic density by
\[
\overline{d}(A) = \limsup_{n \rightarrow \infty}\frac{A(n)}{n}.
\]
The density of a set of natural numbers is defined by
\[
d(A) = \lim_{n \rightarrow \infty}\frac{A(n)}{n},
\]
whenever the limit is exist.

In this paper we prove the existence of a Sidon set such that the lower asymptotic density of its three fold sumset is positive. Namely, we prove the following theorem.

\begin{thm}
There exists a Sidon set $A$ with $\underline{d}(A + A + A) > 0.064$.
\end{thm}

It was proved in \cite{HL}, p.89 that if $A$ is a Sidon set, then 
\[
\liminf_{n\to \infty}\frac{A(n)}{\sqrt{\frac{n}{\log n}}} < \infty, 
\] 
which implies that $\underline{d}(A + A) = 0$. 

It is well known \cite{HL}, p.91-95 that there exists a Sidon set $A$ such that 
\[
\limsup_{n\to \infty}\frac{A(n)}{\sqrt{n}}\ge \frac{\sqrt{2}}{2}, 
\]
which gives $\overline{d}(A + A) \ge \frac{1}{8}$. On the other hand, if $A$ is a Sidon set, then $A(n) \le (1 + o(1))\sqrt{n}$, which implies that
\[
(A + A)(n) \le \left(\frac{1}{2} + o(1)\right)n.
\]
It follows that $\overline{d}(A + A) \le \frac{1}{2}$. 
Now, we have 
\[
\frac{1}{8} \le \sup_{\substack{A \subset \mathbb{N} \\ A \textnormal{ is Sidon }}}\overline{d}(A + A)\le \frac{1}{2}.
\]
\begin{prob}
Do there exist constants $\frac{1}{8} < c_{1} < c_{2} < \frac{1}{2}$ such that 
\[
\displaystyle c_{1} \le \sup_{\substack{A \subset \mathbb{N} \\ A \textnormal{ is Sidon }}}\overline{d}(A + A)\le c_{2}?
\]
\end{prob}
Furthermore, it is natural to study the upper density of the three fold sumset of a Sidon set. The following problem is weaker than Problem 1, but perhaps it is not hopeless to solve.

\begin{prob}
Does there exist a Sidon set $A$ such that $\overline{d}(A + A + A) = 1$? 
\end{prob}

In the next section we give a short summary about the probabilistic method which plays the crucial role in the proof.

\section{Probabilistic tools}

The proof of Theorem 1 is based on the probabilistic method due to Erd\H{o}s and
R\'enyi. There is an excellent summary of this method in the Halberstam - Roth
book \cite{HL}. We use the notation and terminology of this book.
Let $\Omega$ denote the set of strictly increasing sequences of
positive integers. In this paper we denote the probability of an event
$\mathcal{E}$ by $\mathbb{P}(\mathcal{E})$, and the expectation of a random 
variable $\zeta$ by $\mathbb{E}(\zeta)$.

\begin{lem}
Let 
\[
\alpha_{1}, \alpha_{2}, \alpha_{3} \dots{} 
\]
be real numbers satisfying 
\[
0 \le \alpha_{n} \le 1 \hspace*{4mm} (n = 1, 2, \dots{}).
\]
Then there exists a probability space ($\Omega$, $\mathcal{X}$, $\mathbb{P}$) with the following two properties:
\begin{itemize}
\item[(i)] For every natural number $n$, the event $\mathcal{E}^{(n)} = \{\mathcal{A}$:
  $\mathcal{A} \in \Omega$, $n \in \mathcal{A}\}$ is measurable, and
  $\mathbb{P}(\mathcal{E}^{(n)}) = \alpha_{n}$.
\item[(ii)] The events $\mathcal{E}^{(1)}$, $\mathcal{E}^{(2)}$, ... are independent.
\end{itemize}
\end{lem}
See Theorem 13. in \cite{HL}, p. 142. 
We denote the
characteristic function of the event $\mathcal{E}^{(n)}$ by $t_{(\mathcal{A}, n)}$ or we can say the boolean random variable means that:
\[
t_{(\mathcal{A}, n)} = t_{n} =
\left\{
\begin{aligned}
1 \textnormal{, if } n \in \mathcal{A} \\
0 \textnormal{, if } n \notin \mathcal{A}.
\end{aligned} \hspace*{3mm}
\right.
\]


\noindent Furthermore, we denote the number of solutions of
$a_{i_{1}} + a_{i_{2}} + \dots{} + a_{i_{k}} = n$ by $r_{k}(\mathcal{A}, n)$, 
where 
$a_{i_{1}} \in \mathcal{A}$, $a_{i_{2}} \in \mathcal{A}$, ...,$a_{i_{k}} \in \mathcal{A}$, $1 \le a_{i_{1}} < a_{i_{2}} \dots{}
< a_{i_{k}} < n$.  
Thus 
\[
r_{k}(\mathcal{A}, n) = \sum_{\overset{(a_{1}, a_{2}, \dots{}, a_{k}) \in 
\mathbb{N}^{k}}{1
    \le a_{1} < \dots{} < a_{k} < n}\atop {a_{1} + a_{2} + \dots{} + a_{k} =
    n}}t_{(\mathcal{A}, a_{1})}t_{(\mathcal{A}, a_{2})} \dots{} t_{(\mathcal{A}, a_{k})}.
\]
It is easy to see that $r_{k}(\mathcal{A}, n)$ is the sum of random
variables. However for $k > 2$ these variables are not independent because the
same $t_{(\mathcal{A}, a_{i})}$ may appear in many terms. To overcome this
problem we need deeper probabilistic tools.

Our proof is based on a method of J. H. Kim and V. H. Vu. In the next section
we give a short survey of this method. Interested reader can find more
details in \cite{JK}, \cite{TA}, \cite{VA}, \cite{VU}. 
Assume that $t_{1}, t_{2}, \dots{} , t_{n}$ are
independent binary (i.e., all $t_{i}$'s are in $\{0,1\}$) random variables.
Consider a polynomial $Y = Y(t_{1}, t_{2}, \dots{} ,
t_{n})$ in $t_{1}, t_{2}, \dots{} ,
t_{n}$ with degree $k$ (where the degree of this polynomial equals
to the maximum of the sum of the exponents of the monomials.) We say a polynomial $Y$ is totally positive if it can be
written in the form $Y = \sum_{i}e_{i}\Gamma_{i}$, where the $e_{i}$'s are
positive and $\Gamma_{i}$ is a product of some $t_{j}$'s. 
Moreover, $Y$ is regular if all of its coefficients are between zero and one.
We also say $Y$
is simplified, if all of its monomials are square-free
(i.e. do not contain any factor of $t_{i}^{2}$), 
and homogeneous if all the monomials have the same
degree. Thus for instance a boolean polynomial is automatically regular and 
simplified, though not necessarily homogeneous. 
Given any multi-index 
$\underline{\eta} = (\eta_{1}, \dots{} ,\eta_{n}) \in \mathbb{N}^{n}$, we define the partial derivative $\partial^{\underline{\eta}}(Y)$ of
$Y$ by
\[
\partial^{\underline{\eta}}(Y) = \left(\frac{\partial}{\partial t_{1}}\right)^{\eta_{1}} \cdots{} \left(\frac{\partial}{\partial t_{n}}\right)^{\eta_{n}}Y(t_{1}, t_{2}, \dots{} ,t_{n}),
\]
and denote the order of $\underline{\eta}$ as $|\underline{\eta}| = \eta_{1} + \dots{} + \eta_{n}$. For any order $d \ge 0$, we denote $\mathbb{E}_{d}(Y) = \max_{\eta: |\eta| = d}\mathbb{E}(\partial^{\underline{\eta}}(Y))$. Then $\mathbb{E}_{0}(Y) = \mathbb{E}(Y)$ and $\mathbb{E}_{d}(Y) = 0$ if $d$ exceeds the degree of $Y$. 
Define $\mathbb{E}_{\ge d}(Y) =  \max_{d^{'} \ge d}\mathbb{E}_{d^{'}}(Y)$.
The following result was proved by Kim and Vu.

\begin{lem}(J. H. Kim - V. H. Vu)
For every positive integer $k$ and $Y = Y(t_{1}, t_{2}, \dots{} ,t_{n})$ totally 
positive polynomial of degree $k$, where
the $t_{i}$'s are independent binary random variables, and for any 
$\lambda > 0$ there exists a constant $d_{k} > 0$ depending only on $k$ such 
that 
\[
\mathbb{P}\left(|Y-\mathbb{E}(Y)| \ge d_{k}\lambda^{k-\frac{1}{2}}\sqrt{\mathbb{E}_{\ge 0}(Y)\mathbb{E}_{\ge 1}(Y)}\right) = O_{k}\left(e^{-\lambda/4+(k-1)\log n}\right).
\]
\end{lem}

\noindent See \cite{VU} for the proof. Informally, this theorem states that 
when the derivatives of $Y$ are smaller on average than $Y$ itself, and the degree of $Y$ is small, then $Y$ is concentrated
around its mean. Finally we need the Borel - Cantelli lemma:

\begin{lem}(Borel - Cantelli)
Let $\{B_{i}\}$ be a sequence of events in a probability space. If 
\[
\sum_{j=1}^{+\infty}\mathbb{P}(B_{j}) < \infty,
\]
\noindent then with probability 1, at most a finite number of the events
$B_{j}$ can occur.
\end{lem}
This is Theorem 7. in \cite{HL}, p. 135.

\section{Proof of Theorem 4}

\subsection{Outline of the proof}

Let $0 < c < 1$. Define the random sequence $S$ by $\mathbb{P}(n \in S) = \frac{c}{n^{2/3}}$ for every positive integer $n$. Let $T$ be the set
\[
T = \{s \in S: \exists s^{'}, s{''}, s^{'''} < s, s^{'}, s{''}, s^{'''}\in S, s + s^{'} = s{''} + s^{'''}\}.
\]
Then $S \setminus T$ is clearly a Sidon set. On the other hand, 
\[
(S + S + S)(N) \le ((S \setminus T) + (S \setminus T) + (S \setminus T))(N) + (S + S + T)(N), 
\]
thus we have
\[
((S \setminus T) + (S \setminus T) + (S \setminus T))(N) \ge (S + S + S)(N) - (S + S + T)(N). 
\]
It was proved by Goguel \cite{GO} that 
\[
d(S + S + S) = 1 - e^{c^{3}\frac{\Gamma(1/3)^{3}}{6}}, 
\]
with probability 1, where $\Gamma(.)$ is the gamma 
function. We show that with probability 1, 
\begin{equation}
\overline{d}(S + S + T) \le 1.8\Gamma(1/3)^{3}c^{6},
\end{equation}  
Then, with probability 1,
\[
\underline{d}((S \setminus T) + (S \setminus T) + (S \setminus T)) \ge 1 - e^{c^{3}\frac{\Gamma(1/3)^{3}}{6}} - 1.8\Gamma(1/3)^{3}c^{6}.
\]
It is easy to check that $\max_{0 \le c \le 1}1 - e^{c^{3}\frac{\Gamma(1/3)^{3}}{6}} - 1.8\Gamma(1/3)^{3}c^{6} \ge 0.064$, which proves the statement.

To prove (1), we need an upper estimation to $(S + S + T)(N)$. To do this, we verify the following lemma.
\begin{lem}
Let $A,B$ be sequences of positive integers with $A(N) \le (a + o(1))N^{\alpha}$ and $B(N) \le (b + o(1))N^{\beta}$ as $N \rightarrow \infty$. Then 
\[
\sum_{n \le N}R_{A+B}(n) \le \left(\frac{ab\alpha\beta\Gamma(\alpha)\Gamma(\beta)}{(\alpha+\beta)\Gamma(\alpha+\beta)} + o(1)\right)N^{\alpha+\beta} \text{ as } N \rightarrow \infty. 
\]
\end{lem}
Furthermore, it was proved \cite{GO} that with probability 1,
\begin{equation}
S(N) = (3c + o(1))N^{1/3}
\end{equation}
as $N \rightarrow \infty$. 
Next, we show that almost surely
\begin{equation}
T(N) \le (10.8c^{4} + o(1))N^{1/3}
\end{equation}
as $N \rightarrow \infty$. If $n \in S + S$ and $\frac{n}{2} \notin S$, then $R_{S+S}(n) \ge 2$. By Lemma 4 and (2) we have
\[
(S + S)(N) \le \frac{1}{2}\left(\sum_{n \le N}R_{S+S}(n) + S\left(\frac{N}{2}\right)\right) \le \left(\frac{1}{2}\cdot \frac{(3c)(3c)\cdot \frac{1}{3} \cdot \frac{1}{3} \Gamma(\frac{1}{3})\Gamma(\frac{1}{3})}
{(\frac{1}{3}+\frac{1}{3})\Gamma(\frac{1}{3}+\frac{1}{3})} + o(1)\right)N^{\frac{2}{3}} 
\]
\[
= \left(\frac{3}{4}\cdot \frac{\Gamma(\frac{1}{3})^{2}}{\Gamma(\frac{2}{3})}c^{2} + o(1)\right)N^{\frac{2}{3}}
\]
as $N \rightarrow \infty$. By using Lemma 4 and (3), one get that with probability 1, 
\[
(S + S + T)(N) \le \left(\frac{\frac{3}{4}\cdot \frac{\Gamma(\frac{1}{3})^{2}}{\Gamma(\frac{2}{3})}\cdot c^{2}\cdot 10.8c^{4}\cdot \frac{1}{3} \cdot \frac{2}{3} \Gamma(\frac{1}{3})\Gamma(\frac{2}{3})}{(\frac{1}{3}+\frac{2}{3})\Gamma(\frac{1}{3}+\frac{2}{3})} + o(1) \right)N
\] 
\[
= \left(1.8\Gamma\left(\frac{1}{3}\right)^{3}c^{6} + o(1)\right)N
\]
as $N \rightarrow \infty$, which proves (1). 

\subsection{Proof of Lemma 4}

Let $\varepsilon, \varepsilon_{1}, \varepsilon_{2}, \varepsilon_{3} > 0$ and $A = \{a_{1}, a_{2}, \dots{},\}$, $(1 \le a_{1} < a_{2} < \dots{})$. Then $n = A(a_{n}) \le (a + o(1))a_{n}^{\alpha}$, thus
$a_{n} \ge \left(\frac{1}{a^{1/\alpha}} + o(1)\right)n^{1/\alpha}$ as $n \rightarrow \infty$. Then there exists an $N_{1}(\varepsilon_{1})$ such that if $n \ge N_{1}(\varepsilon_{1})$, then
$a_{n} \ge \left(\frac{1}{a^{1/\alpha}} - \varepsilon_{1} \right)n^{1/\alpha}$. It follows that 
\[
\sum_{n \le N}R_{A+B}(n) = \sum_{a_{n} < N}B(N-a_{n}) \le N_{1}(\varepsilon_{1})B(N) + \sum_{n: \left(\frac{1}{a^{1/\alpha}} - \varepsilon_{1} \right)n^{1/\alpha} < N}B\left(N - \left(\frac{1}{a^{1/\alpha}} - \varepsilon_{1} \right)n^{1/\alpha}\right).
\]
There exists an $N_{2}(\varepsilon_{2})$ such that if $n \ge N_{2}(\varepsilon_{2})$, then $B(n) < (b + \varepsilon_{2})n^{\beta}$.
This implies that
\[
\sum_{n \le N}R_{A+B}(n) \le O(B(N)) + (b + \varepsilon_{2})\sum_{n: \left(\frac{1}{a^{1/\alpha}} - \varepsilon_{1} \right)n^{1/\alpha} < N}\left(N - \left(\frac{1}{a^{1/\alpha}} - \varepsilon_{1} \right)n^{1/\alpha}\right)^{\beta} 
\]
\[
= O(B(N)) + (b + \varepsilon_{2})N^{\beta}\cdot \frac{N^{\alpha}}{\left(\frac{1}{a^{1/\alpha}} - \varepsilon_{1} \right)^{\alpha}}\cdot \frac{1}{\frac{N^{\alpha}}{\left(\frac{1}{a^{1/\alpha}} - \varepsilon_{1} \right)^{\alpha}}}\sum_{n \le \frac{N^{\alpha}}{\left(\frac{1}{a^{1/\alpha}} - \varepsilon_{1} \right)^{\alpha}}}\left(1 - \left(\frac{n}{\frac{N^{\alpha}}{\left(\frac{1}{a^{1/\alpha}} - \varepsilon_{1} \right)^{\alpha}}}\right)^{1/\alpha}\right)^{\beta}.
\]
Moreover, it is clear that there exists an $N_{3}(\varepsilon_{3})$ such that if $n \ge N_{3}(\varepsilon_{3})$, then 
\[
\frac{1}{\frac{N^{\alpha}}{\left(\frac{1}{a^{1/\alpha}} - \varepsilon_{1} \right)^{\alpha}}}\sum_{n \le \frac{N^{\alpha}}{\left(\frac{1}{a^{1/\alpha}} - \varepsilon_{1} \right)^{\alpha}}}\left(1 - \left(\frac{n}{\frac{N^{\alpha}}{\left(\frac{1}{a^{1/\alpha}} - \varepsilon_{1} \right)^{\alpha}}}\right)^{1/\alpha}\right)^{\beta} \le \int_{0}^{1}(1 - x^{1/\alpha})^{\beta}dx + \varepsilon_{3}.
\]
Now, we compute $\int_{0}^{1}(1 - x^{1/\alpha})^{\beta}dx$. Writing $y = x^{1/\alpha}$, we have
\[
\int_{0}^{1}(1 - x^{1/\alpha})^{\beta}dx = \int_{0}^{1}(1 - y)^{\beta}\alpha\cdot y^{\alpha-1}dy = \alpha\cdot \frac{\Gamma(\alpha)\Gamma(\beta+1)}{\Gamma(\alpha+\beta+1)} = \frac{\alpha\beta\Gamma(\alpha)\Gamma(\beta)}{(\alpha+\beta)\Gamma(\alpha+\beta)},  
\]  
where $\Gamma(x)$ is the Gamma function. This implies that
\[
\sum_{n \le N}R_{A+B}(n) \le O(B(N)) + \frac{1}{(a^{1/\alpha} - \varepsilon_{1})^{\alpha}}(b + \varepsilon_{2})\left(\frac{\alpha\beta\Gamma(\alpha)\Gamma(\beta)}{(\alpha+\beta)\Gamma(\alpha+\beta)} + \varepsilon_{3} \right)N^{\alpha+\beta}.
\]
By choosing $\varepsilon_{1}, \varepsilon_{2}, \varepsilon_{3}$ small enough depending on $\varepsilon$, then there exists an $N_{0}(\varepsilon)$ such that if $n \ge N_{0}(\varepsilon)$, then
\[
\sum_{n \le N}R_{A+B}(n) \le \left(\frac{ab\alpha\beta\Gamma(\alpha)\Gamma(\beta)}{(\alpha+\beta)\Gamma(\alpha+\beta)} + \varepsilon \right)N^{\alpha+\beta}.
\]
The proof of the lemma is completed.

\subsection{Upper estimation for $T(N)$}

Now we prove (3). We write
\[
t_{n} = 
\left\{
\begin{aligned}
1 \textnormal{, if } n \in S \\
0 \textnormal{, if } n \notin S.
\end{aligned} \hspace*{3mm}
\right.
\]
Let
\[
Z(N) = \sum_{\overset{(x_{1}, x_{2}, x_{3}, x_{4}) \in \mathbb{N}^{4}}{1 \le x_{4} 
< x_{3} < x_{2} < x_{1} \le N}\atop {x_{1} + x_{4} = x_{2} + x_{3}}}t_{x_{1}}t_{x_{2}}t_{x_{3}}t_{x_{4}} + \sum_{\overset{(x_{1}, x_{2}, x_{4}) \in \mathbb{N}^{4}}{1 \le x_{4} < x_{2} < x_{1} \le N}\atop {x_{1} + x_{4} = 2x_{2}}}t_{x_{1}}t_{x_{2}}t_{x_{4}}.
\]
It is clear that $T(N) \le Z(N)$, thus it is enough to prove that with probability 1, 
\begin{equation}
Z(N) \le (10.8c^{4} + o(1))N^{1/3}.
\end{equation}
as $N \rightarrow \infty$. To prove (4), we show that 
\begin{equation}
\mathbb{E}(Z(N)) \le (10.8c^{4} + o(1))N^{1/3} 
\end{equation}
as $N \rightarrow \infty$ and with probability 1, 
\begin{equation}
Z(N) = \mathbb{E}(Z(N)) + O(N^{1/6}\log^{4}N).
\end{equation}
First, we prove (5). Clearly we have
\[
\mathbb{E}(Z(N)) = \mathbb{E}\left(\sum_{\substack{(x_{1}, x_{2}, x_{4}) \in \mathbb{N}^{4}\\1 \le x_{4} <  x_{2} < x_{1} \le N \\x_{1} + x_{4} = 2x_{2}}}t_{x_{1}}t_{x_{2}}t_{x_{4}} \right) + \mathbb{E}\left(\sum_{\overset{(x_{1}, x_{2}, x_{3}, x_{4}) \in \mathbb{N}^{4}}{1 \le x_{4} 
< x_{3} < x_{2} < x_{1} \le N}\atop {x_{1} + x_{4} = x_{2} + x_{3}}}t_{x_{1}}t_{x_{2}}t_{x_{3}}t_{x_{4}} \right) = \mathbb{E}_{1} + \mathbb{E}_{2}.
\]
Then we have 
\[
\mathbb{E}_{1} = \sum_{x_{1} = 3}^{N}\frac{c}{x_{1}^{2/3}}\sum_{\overset{1 \le x_{4} < x_{1}}{x_{1} \equiv x_{4} \bmod{2}}}\frac{c}{x_{4}^{2/3}}\cdot \frac{c}{\left(\frac{x_{1}+x_{4}}{2}\right)^{2/3}}.
\]
Since
\[
\frac{1}{\left(\frac{x_{1}+x_{4}}{2}\right)^{2/3}} \le \frac{1}{\left(\frac{x_{1}}{2}\right)^{2/3}} = \frac{2^{2/3}}{x_{1}^{2/3}},
\]
and by the Euler formula for the estimation of a sum by integral, we have
\[
\mathbb{E}_{1} = O\left(\sum_{x_{1} = 1}^{N}\frac{1}{x_{1}^{2/3}}\sum_{1 \le x_{4} \le x_{1}}\frac{1}{x_{4}^{2/3}}\cdot \frac{1}{x_{1}^{2/3}}\right) = O\left(\sum_{x_{1} = 1}^{N}\frac{1}{x_{1}^{4/3}}\sum_{1 \le x_{4} \le x_{1}}\frac{1}{x_{4}^{2/3}}\right) 
\]
\[
= O\left(\sum_{x_{1} = 1}^{N}\frac{1}{x_{1}}\right) = O(\log N).
\]
On the other hand,
\[
\mathbb{E}_{2} = \sum_{x_{1} = 2}^{N}\frac{c}{x_{1}^{2/3}}\sum_{x_{4} = 1}^{x_{1}-1}\frac{c}{x_{4}^{2/3}}\sum_{x_{4} < x_{2} < \frac{x_{1}+x_{4}}{2}}\frac{c}{x_{2}^{2/3}}\cdot \frac{c}{(x_{1}+x_{4}-x_{2})^{2/3}}
\]
\[
= c^{4}N^{1/3}\cdot \frac{1}{N^{3}}\sum_{x_{1} = 2}^{N}\frac{1}{\left(\frac{x_{1}}{N}\right)^{2/3}}\sum_{x_{4} = 1}^{x_{1}-1}\frac{1}{\left(\frac{x_{4}}{N}\right)^{2/3}}\sum_{x_{4} < x_{2} < \frac{x_{1}+x_{4}}{2}}\frac{1}{\left(\frac{x_{2}}{N}\right)^{2/3}}\cdot \frac{1}{\left(\frac{x_{1}}{N}+\frac{x_{4}}{N}-\frac{x_{2}}{N}\right)^{2/3}}.
\]
Consider the open set $D \subseteq \mathbb{R}^{3}$
\[
D = \left\{(x,y,z): 0 < x < 1, 0 < y < x, y < z < \frac{x+y}{2}\right\}
\]
and the unbounded function $f: D \rightarrow \mathbb{R}$,
\[
f(x,y,z) = \frac{1}{x^{2/3}}\cdot \frac{1}{y^{2/3}}\cdot \frac{1}{z^{2/3}}\cdot \frac{1}{(x+y-z)^{2/3}}.
\]
Now we show that the integral $\iiint_{D}f(x,y,z)dV$ exists.
Let $M > 0$, $f_{M}(x,y,z) = \min\{f(x,y,z),M\}$ if $(x,y,z) \in D$ and let
$I(M) = \iiint_{D}f_{M}(x,y,z)dV$. We have to prove that $\lim_{M\rightarrow \infty}\iiint_{D}f_{M}(x,y,z)dV$ exists. Since $I(M)$ is monotone increasing, it is enough to show that $I(M)$ is bounded above. Let
\[
D_{M} = \left\{(x,y,z): 0 < x < 1, 0 < y < x, y < z < \frac{x+y}{2}, y < \frac{1}{M^{2}}\right\}.
\]
For $(x,y,z) \in D_{M}$, we have $0 < x < 1$, $0 < y < \frac{1}{M^{2}}$, 
$0 < z < \frac{x+y}{2} < 1$, thus the volume of $D_{M}$ is $\text{Vol}(D_{M}) \le \frac{1}{M^{2}}$ and so
$\iiint_{D_{M}}f_{M}(x,y,z)dV \le M\cdot \frac{1}{M^{2}} = \frac{1}{M}$. Moreover,
\[
D \setminus D_{M} = \left\{(x,y,z): 0 < x < 1, 0 < y < x, y < z < \frac{x+y}{2}, y \ge \frac{1}{M^{2}}\right\}
\]
\[
= \left\{(x,y,z): \frac{1}{M^{2}} < x < 1, \frac{1}{M^{2}} \le y < x, y < z < \frac{x+y}{2}\right\},
\]
then
\[
\iiint_{D\setminus D_{M}}f_{M}(x,y,z)dV \le \int_{\frac{1}{M^{2}}}^{1}\left(\int_{\frac{1}{M^{2}}}^{x}\left(\int_{y}^{\frac{x+y}{2}}\frac{1}{x^{2/3}}\cdot \frac{1}{y^{2/3}}\cdot \frac{1}{z^{2/3}}\cdot \frac{1}{(x+y-z)^{2/3}}dz\right)dy\right)dx.
\]
Since $z < \frac{x+y}{2}$ and $0 < y < x$, thus $\frac{1}{(x+y-z)^{2/3}} \le \frac{1}{\left(\frac{x+y}{2}\right)^{2/3}} = \frac{2^{2/3}}{(x+y)^{2/3}} < \frac{2^{2/3}}{x^{2/3}}$, and so
\[
\iiint_{D\setminus D_{M}}f_{M}(x,y,z)dV \le \int_{\frac{1}{M^{2}}}^{1}\frac{2^{2/3}}{x^{4/3}}\left(\int_{\frac{1}{M^{2}}}^{x}\frac{1}{y^{2/3}}\left(\int_{y}^{x}\frac{1}{z^{2/3}}dz\right)dy\right)dx
\]
\[
\le \int_{\frac{1}{M^{2}}}^{1}\frac{3\cdot 2^{2/3}}{x}\left(\int_{\frac{1}{M^{2}}}^{x}\frac{1}{y^{2/3}}dy\right)dx \le \int_{\frac{1}{M^{2}}}^{1}\frac{9\cdot 2^{2/3}}{x^{2/3}}dx < 27\cdot 2^{2/3}. 
\]
This implies that $I(M) < 27\cdot 2^{2/3} + \frac{1}{M}$, and so $I(M) \le 27\cdot 2^{2/3}$ because $I(M)$ is monotone increasing.

Next, we prove that 
\begin{equation}
\lim_{N \rightarrow \infty}\frac{1}{N^{3}}\sum_{x_{1} = 2}^{N}\frac{1}{\left(\frac{x_{1}}{N}\right)^{2/3}}\sum_{x_{4} = 1}^{x_{1}-1}\frac{1}{\left(\frac{x_{4}}{N}\right)^{2/3}}\sum_{x_{4} < x_{2} < \frac{x_{1}+x_{4}}{2}}\frac{1}{\left(\frac{x_{2}}{N}\right)^{2/3}}\cdot \frac{1}{\left(\frac{x_{1}}{N}+\frac{x_{4}}{N}-\frac{x_{2}}{N}\right)^{2/3}} 
\end{equation}
\[
= \iiint_{D}f(x,y,z)dV.
\]

It is clear that 
\[
\frac{1}{N^{3}}\sum_{x_{1} = 2}^{N}\sum_{x_{4} = 1}^{x_{1}-1}\sum_{x_{1} < x_{2} < \frac{x_{1}+x_{4}}{2}}\frac{1}{\left(\frac{x_{1}}{N}\right)^{2/3}}\cdot \frac{1}{\left(\frac{x_{4}}{N}\right)^{2/3}}\cdot \frac{1}{\left(\frac{x_{1}}{N}+\frac{x_{4}}{N}-\frac{x_{2}}{N}\right)^{2/3}} 
\]
\[
= \frac{1}{N^{3}}\sum_{x_{1} = 2}^{N}\sum_{x_{4} = 1}^{x_{1}-1}\sum_{x_{1} < x_{2} < \frac{x_{1}+x_{4}}{2}}f\left(\frac{x_{1}}{N},\frac{x_{4}}{N},\frac{x_{2}}{N}\right).
\]
Let $\varepsilon > 0$. Clearly we have
\[
\left|\frac{1}{N^{3}}\sum_{x_{1} = 2}^{N}\sum_{x_{4} = 1}^{x_{1}-1}\sum_{x_{4} < x_{2} < \frac{x_{1}+x_{4}}{2}}f\left(\frac{x_{1}}{N},\frac{x_{4}}{N},\frac{x_{2}}{N}\right) - \iiint_{D}f(x,y,z)dV\right| 
\]
\[
\le
\left|\iiint_{D}f(x,y,z)dV - \iiint_{D}f_{M}(x,y,z)dV\right| 
\]
\[
+ \left|\iiint_{D}f_{M}(x,y,z)dV - \frac{1}{N^{3}}\sum_{x_{1} = 2}^{N}\sum_{x_{4} = 1}^{x_{1}-1}\sum_{x_{4} < x_{2} < \frac{x_{1}+x_{4}}{2}}f_{M}\left(\frac{x_{1}}{N},\frac{x_{4}}{N},\frac{x_{2}}{N}\right)\right| 
\]
\[
+ \left|\frac{1}{N^{3}}\sum_{x_{1} = 2}^{N}\sum_{x_{4} = 1}^{x_{1}-1}\sum_{x_{4} < x_{2} < \frac{x_{1}+x_{4}}{2}}f_{M}\left(\frac{x_{1}}{N},\frac{x_{4}}{N},\frac{x_{2}}{N}\right) - \frac{1}{N^{3}}\sum_{x_{1} = 2}^{N}\sum_{x_{4} = 1}^{x_{1}-1}\sum_{x_{4} < x_{2} < \frac{x_{1}+x_{4}}{2}}f\left(\frac{x_{1}}{N},\frac{x_{4}}{N},\frac{x_{2}}{N}\right)\right|.
\]
Obviously,
\[
\left|\frac{1}{N^{3}}\sum_{x_{1} = 2}^{N}\sum_{x_{4} = 1}^{x_{1}-1}\sum_{x_{4} < x_{2} < \frac{x_{1}+x_{4}}{2}}f_{M}\left(\frac{x_{1}}{N},\frac{x_{4}}{N},\frac{x_{2}}{N}\right) - \frac{1}{N^{3}}\sum_{x_{1} = 2}^{N}\sum_{x_{4} = 1}^{x_{1}-1}\sum_{x_{4} < x_{2} < \frac{x_{1}+x_{4}}{2}}f\left(\frac{x_{1}}{N},\frac{x_{4}}{N},\frac{x_{2}}{N}\right)\right|
\]
\[
= \frac{1}{N^{3}}\sum_{\substack{(x_{1},x_{2},x_{4})\\1 \le x_{1} \le N\\1 \le x_{4} < x_{1}\\x_{4} < x_{2} < \frac{x_{1}+x_{4}}{2}\\f\left(\frac{x_{1}}{N},\frac{x_{4}}{N},\frac{x_{2}}{N}\right)>M}}\left(f\left(\frac{x_{1}}{N},\frac{x_{4}}{N},\frac{x_{2}}{N}\right) - M\right) \le \frac{1}{N^{3}}\sum_{\substack{(x_{1},x_{2},x_{4})\\1 \le x_{1} \le N\\1 \le x_{4} < x_{1}\\x_{4} < x_{2} < \frac{x_{1}+x_{4}}{2}\\f\left(\frac{x_{1}}{N},\frac{x_{4}}{N},\frac{x_{2}}{N}\right) > M}}f
\left(\frac{x_{1}}{N},\frac{x_{4}}{N},\frac{x_{2}}{N}\right).
\]
It is clear that 
\[
f\left(\frac{x_{1}}{N},\frac{x_{4}}{N},\frac{x_{2}}{N}\right) > M
\]
if and only if 
\[
\frac{1}{\left(\frac{x_{1}}{N}\right)^{2/3}}\cdot \frac{1}{\left(\frac{x_{4}}{N}\right)^{2/3}}\cdot \frac{1}{\left(\frac{x_{2}}{N}\right)^{2/3}}\cdot \frac{1}{\left(\frac{x_{1}}{N}+\frac{x_{4}}{N}-\frac{x_{2}}{N}\right)^{2/3}} > M
\]
if and only if $\frac{N^{4}}{M^{3/2}} > x_{1}x_{4}x_{2}(x_{1} + x_{4} - x_{2})$. Since $x_{4} < x_{1}$, $x_{1} + x_{4} - x_{2} > \frac{x_{1} + x_{4}}{2} > x_{4}$, then  
$\frac{N^{4}}{M^{3/2}} > x_{4}^{4}$, i.e., $\frac{N}{M^{3/8}} > x_{4}$. 
Then we have
\[
\frac{1}{N^{3}}\sum_{\substack{(x_{1},x_{4},x_{2})\\1 \le x_{1} \le N\\1 \le x_{4} < x_{1}\\x_{4} < x_{2} < \frac{x_{1}+x_{4}}{2}\\f\left(\frac{x_{1}}{N},\frac{x_{4}}{N},\frac{x_{2}}{N}\right)>M}}f\left(\frac{x_{1}}{N},\frac{x_{4}}{N},\frac{x_{2}}{N}\right)    
\le \frac{1}{N^{3}}\sum_{\substack{(x_{1},x_{4},x_{2})\\1 \le x_{1} \le N\\1 \le x_{4} \le \min\left\{x_{1},\frac{N}{M^{3/8}}\right\}\\x_{4} < x_{2} < \frac{x_{1}+x_{4}}{2}}}f\left(\frac{x_{1}}{N},\frac{x_{4}}{N},\frac{x_{2}}{N}\right)
\]
\[
= \frac{1}{N^{3}}\sum_{x_{1} = 1}^{N}\sum_{x_{4} = 1}^{\min\left\{x_{1},\frac{N}{M^{3/8}}\right\}}\sum_{x_{4} < x_{2} < \frac{x_{1}+x_{4}}{2}}\frac{1}{\left(\frac{x_{1}}{N}\right)^{2/3}}\cdot \frac{1}{\left(\frac{x_{4}}{N}\right)^{2/3}}\cdot \frac{1}{\left(\frac{x_{2}}{N}\right)^{2/3}}\cdot \frac{1}{\left(\frac{x_{1}}{N}+\frac{x_{4}}{N}-\frac{x_{2}}{N}\right)^{2/3}} 
\]
\[
= \frac{1}{N^{1/3}}\sum_{x_{1} = 1}^{N}\frac{1}{x_{1}^{2/3}}\sum_{x_{4} = 1}^{\min\left\{x_{1},\frac{N}{M^{3/8}}\right\}}\frac{1}{x_{4}^{2/3}}\sum_{x_{4} < x_{2} < \frac{x_{1}+x_{4}}{2}}\frac{1}{x_{2}^{2/3}}\frac{1}{(x_{1}+x_{4}-x_{2})^{2/3}}.
\]
Obviously, 
\[
\frac{1}{(x_{1}+x_{4}-x_{2})^{2/3}} < \frac{1}{\left(\frac{x_{1}+x_{4}}{2}\right)^{2/3}} = \frac{2^{2/3}}{(x_{1}+x_{4})^{2/3}} < \frac{2^{2/3}}{x_{1}^{2/3}}, 
\]
thus we have 
\[
\frac{1}{N^{3}}\sum_{\substack{(x_{1},x_{2},x_{4})\\1 \le x_{1} \le N\\1 \le x_{4} < x_{1}\\x_{4} < x_{2} < \frac{x_{1}+x_{4}}{2}\\f\left(\frac{x_{1}}{N},\frac{x_{4}}{N},\frac{x_{2}}{N}\right)>M}}f\left(\frac{x_{1}}{N},\frac{x_{4}}{N},\frac{x_{2}}{N}\right) 
< \frac{1}{N^{1/3}}\sum_{x_{1} = 1}^{N}\frac{2^{2/3}}{x_{1}^{4/3}}\sum_{x_{4} = 1}^{\min\left\{x_{1},\frac{N}{M^{3/8}}\right\}}\frac{1}{x_{4}^{2/3}}\sum_{x_{2} = 1}^{x_{1}}\frac{1}{x_{2}^{2/3}} 
\]
\[
< \frac{1}{N^{1/3}}\sum_{x_{1} = 1}^{N}\frac{2^{2/3}}{x_{1}^{4/3}}\sum_{x_{4} = 1}^{\min\left\{x_{1},\frac{N}{M^{3/8}}\right\}}\frac{1}{x_{4}^{2/3}}\cdot 3x_{1}^{1/3} = \frac{3\cdot 2^{2/3}}{N^{1/3}}\sum_{x_{1} = 1}^{\left\lfloor \frac{N}{M^{3/8}}\right\rfloor}\frac{1}{x_{1}}
\sum_{x_{4} = 1}^{x_{1}}\frac{1}{x_{4}^{2/3}} 
\]
\[
+ \frac{3\cdot 2^{2/3}}{N^{1/3}}\sum_{x_{1} =\left\lfloor \frac{N}{M^{3/8}}\right\rfloor+1}^{N}\frac{1}{x_{1}}\sum_{x_{4} = 1}^{\left\lfloor \frac{N}{M^{3/8}}\right\rfloor}\frac{1}{x_{4}^{2/3}}
\]
\[
< \frac{3\cdot 2^{2/3}}{N^{1/3}}\sum_{x_{1} = 1}^{\left\lfloor \frac{N}{M^{3/8}}\right\rfloor}\frac{1}{x_{1}}\cdot 3x_{1}^{1/3} + \frac{3\cdot 2^{2/3}}{N^{1/3}}\sum_{x_{1} =\left\lfloor \frac{N}{M^{3/8}}\right\rfloor+1}^{N}\frac{1}{x_{1}}\cdot 3\cdot \left(\frac{N}{M^{3/8}}\right)^{1/3}
\]
\[
< \frac{9\cdot 2^{2/3}}{N^{1/3}}\sum_{x_{1} = 1}^{\left\lfloor \frac{N}{M^{3/8}}\right\rfloor}\frac{1}{x_{1}^{2/3}} + \frac{9\cdot 2^{2/3}}{M^{1/8}}\sum_{x_{1} =\left\lfloor\frac{N}{M^{3/8}}\right\rfloor+1}^{N}\frac{1}{x_{1}} 
< \frac{9\cdot 2^{2/3}}{N^{1/3}}\cdot 3\cdot \left(\frac{N}{M^{3/8}}\right)^{1/3} 
\]
\[
+ \frac{9\cdot 2^{2/3}}{M^{1/8}}\left(\log N + 1 - \log\frac{N}{M^{3/8}}
\right) 
= \frac{27\cdot 2^{2/3}}{M^{1/8}} + \frac{9\cdot 2^{2/3}(1+\frac{3}{8}\log M)}{M^{1/8}}. 
\]
Then there exists an $M_{0}(\varepsilon) > 0$ such that if $M \ge M_{0}(\varepsilon)$, then
\[
\left|\frac{1}{N^{3}}\sum_{x_{1} = 1}^{N}\sum_{x_{4} = 1}^{x_{1}-1}\sum_{x_{4} < x_{2} < \frac{x_{1}+x_{4}}{2}}f_{M}\left(\frac{x_{1}}{N},\frac{x_{4}}{N},\frac{x_{2}}{N}\right) - \frac{1}{N^{3}}\sum_{x_{1} = 1}^{N}\sum_{x_{4} = 1}^{x_{2}-1}\sum_{x_{4} < x_{2} < \frac{x_{1}+x_{4}}{2}}f\left(\frac{x_{1}}{N},\frac{x_{4}}{N},\frac{x_{2}}{N}\right)\right| < \frac{\varepsilon}{3}.
\]
Since $\lim_{M \rightarrow \infty}\iiint_{D}f_{M}(x,y,z)dV = \iiint_{D}f(x,y,z)dV$, then we choose such an $M > M_{0}(\varepsilon)$
\[
\left|\iiint_{D}f(x,y,z)dV - \iiint_{D}f_{M}(x,y,z)dV\right| < \frac{\varepsilon}{3}.
\]
It is clear that there exists an $N_{0}(\varepsilon) > 0$ such that if $N \ge N_{0}(\varepsilon)$, then
\[
\left|\iiint_{D}f_{M}(x,y,z)dV - \frac{1}{N^{3}}\sum_{x_{1} = 2}^{N}\sum_{x_{4} = 1}^{x_{1}-1}\sum_{x_{4} < x_{2} < \frac{x_{1}+x_{4}}{2}}f_{M}\left(\frac{x_{1}}{N},\frac{x_{4}}{N},\frac{x_{2}}{N}\right)\right| < \frac{\varepsilon}{3}. 
\]
Then for $N \ge N_{0}(\varepsilon)$, we have
\[
\left|\frac{1}{N^{3}}\sum_{x_{1} = 1}^{N}\frac{1}{\left(\frac{x_{1}}{N}\right)^{2/3}}\sum_{x_{4} = 1}^{x_{1} - 1}\frac{1}{\left(\frac{x_{4}}{N}\right)^{2/3}}\sum_{x_{4} < x_{2} < \frac{x_{1}+x_{4}}{2}}\frac{1}{\left(\frac{x_{2}}{N}\right)^{2/3}}\cdot \frac{1}{\left(\frac{x_{1}}{N}+\frac{x_{4}}{N}-\frac{x_{2}}{N}\right)^{2/3}}-
\iiint_{D}f(x,y,z)dV\right| < \varepsilon,
\]
which proves (7). 

Then, we get that $\mathbb{E}_{2} = (c^{4}\iiint_{D}f(x,y,z)dV + o(1))N^{1/3}$ as $N \rightarrow \infty$. By using Matlab, we get that $10.7 < \iiint_{D}f(x,y,z)dV < 10.8$, which proves (5).


To prove (6), we apply the Kim-Vu inequality. We show that there exist constants $c_{1}$ and $c_{2}$ that

\begin{equation}
\mathbb{E}_{\ge 0}(Z(N)) \le c_{1}N^{1/3}
\end{equation}
and
\begin{equation}
\mathbb{E}_{\ge 1}(Z(N)) \le c_{2}\log N.
\end{equation}
To do this, we give upper estimations for the expectation of the partial derivatives of $Z(N)$. Let $1 \le x_{0} \le x$. We will estimate the following sums by integrals. We have
\[
\mathbb{E}\left(\frac{\partial}{\partial t_{x_{0}}} Z(N)\right) = \sum_{\substack{(x_{2},x_{3},x_{4})\\x_{0}+x_{2}=x_{3}+x_{4}\\1 \le x_{2} < x_{0}\\x_{3} < x_{4}}}\frac{c}{x_{2}^{2/3}}\cdot \frac{c}{x_{3}^{2/3}}\cdot \frac{c}{x_{4}^{2/3}} + \sum_{\substack{(x_{2},x_{3},x_{4})\\x_{0}+x_{2}=x_{3}+x_{4}\\x_{0} < x_{2} \le N\\x_{3} < x_{4}}}\frac{c}{x_{2}^{2/3}}\cdot \frac{c}{x_{3}^{2/3}}\cdot \frac{c}{x_{4}^{2/3}} 
\]
\[
+ \sum_{\substack{(x_{3},x_{4})\\2x_{0}=x_{3}+x_{4}\\x_{3} < x_{4}}}\frac{c}{x_{3}^{2/3}}\cdot \frac{c}{x_{4}^{2/3}} +
\sum_{\substack{(x_{2},x_{3})\\x_{0}+x_{2}=2x_{3}\\1 \le x_{2} < x_{0}\\x_{0} \equiv x_{2} \bmod{2}}}\frac{c}{x_{2}^{2/3}}\cdot \frac{c}{x_{3}^{2/3}} 
+ \sum_{\substack{(x_{2},x_{3})\\x_{0}+x_{2}=2x_{3}\\x_{0} < x_{2} \le N\\x_{0} \equiv x_{2} \bmod{2}}}\frac{c}{x_{2}^{2/3}}\cdot \frac{c}{x_{3}^{2/3}}
\]
\[
= \mathbb{E}_{3} + \mathbb{E}_{4} + \mathbb{E}_{5} + \mathbb{E}_{6} + \mathbb{E}_{7}.
\]
In $\mathbb{E}_{3}$, we have $x_{4} > \frac{x_{0}+x_{2}}{2} > \frac{x_{0}}{2}$, 
and so $\frac{1}{x_{4}^{2/3}} < \frac{1}{\left(\frac{x_{0}}{2}\right)^{2/3}} = \frac{2^{2/3}}{x_{0}^{2/3}}$, which implies that
\[
\mathbb{E}_{3} = O\left(\sum_{x_{2}=1}^{x_{0}}\frac{1}{x_{2}^{2/3}}\left(\sum_{x_{3}=1}^{x_{0}}\frac{1}{x_{3}^{2/3}}\cdot \frac{1}{x_{0}^{2/3}}\right)\right) = O\left(\frac{1}{x_{0}^{2/3}}\cdot \sum_{x_{2}=1}^{x_{0}}\frac{1}{x_{2}^{2/3}}\left(\sum_{x_{3}=1}^{x_{0}}\frac{1}{x_{3}^{2/3}}\right)\right) 
\]
\[
= O\left(\frac{1}{x_{0}^{1/3}}\cdot \sum_{x_{2}=1}^{x_{0}}\frac{1}{x_{2}^{2/3}}\right) = O(1).
\]
In $\mathbb{E}_{4}$, we have $x_{4} > \frac{x_{0}+x_{2}}{2} > \frac{x_{2}}{2}$, and so $\frac{1}{x_{4}^{2/3}} < \frac{1}{\left(\frac{x_{2}}{2}\right)^{2/3}} = \frac{2^{2/3}}{x_{2}^{2/3}}$, which implies that
\[
\mathbb{E}_{4} = O\left(\sum_{x_{2}=x_{0}+1}^{N}\frac{1}{x_{2}^{2/3}}\left(\sum_{x_{3}=1}^{x_{2}}\frac{1}{x_{3}^{2/3}}\cdot \frac{1}{x_{2}^{2/3}}\right)\right) = O\left(\sum_{x_{2}=x_{0}+1}^{N}\frac{1}{x_{2}^{4/3}}\left(\sum_{x_{3}=1}^{x_{2}}\frac{1}{x_{3}^{2/3}}\right)\right) 
\]
\[
= O\left(\sum_{x_{2}=x_{0}+1}^{N}\frac{1}{x_{2}}\right) = O(\log N).
\]
In $\mathbb{E}_{5}$, we have $x_{4} > x_{0}$, and so $\frac{1}{x_{4}^{2/3}} < \frac{1}{x_{0}^{2/3}}$, which implies that
\[
\mathbb{E}_{5} = O\left(\sum_{x_{3}=1}^{x_{0}}\frac{1}{x_{3}^{2/3}}\cdot \frac{1}{x_{0}^{2/3}}\right) = O\left(\frac{1}{x_{0}^{1/3}}\right) = O(1). 
\]
In $\mathbb{E}_{6}$, $x_{3} = \frac{x_{0}+x_{2}}{2} > \frac{x_{0}}{2}$, thus we have $\frac{1}{x_{3}^{2/3}} \le \frac{1}{\left(\frac{x_{0}}{2}\right)^{2/3}} = \frac{2^{2/3}}{x_{0}^{2/3}}$, which implies that
\[
\mathbb{E}_{6} = O\left(\sum_{x_{2}=1}^{x_{0}}\frac{1}{x_{2}^{2/3}}\cdot \frac{1}{x_{0}^{2/3}}\right) = O\left(\frac{1}{x_{0}^{1/3}}\right) = O(1). 
\]
In $\mathbb{E}_{7}$, $x_{3} = \frac{x_{0}+x_{2}}{2} > \frac{x_{2}}{2}$, thus we have $\frac{1}{x_{3}^{2/3}} \le \frac{1}{\left(\frac{x_{2}}{2}\right)^{2/3}} = \frac{2^{2/3}}{x_{2}^{2/3}}$, which implies that
\[
\mathbb{E}_{7} = O\left(\sum_{x_{2}=x_{0}+1}^{N}\frac{1}{x_{2}^{2/3}}\cdot \frac{1}{x_{2}^{2/3}}\right) = O(1). 
\]
In summary, we have 
\begin{equation}
\mathbb{E}\left(\frac{\partial}{\partial t_{x_{0}}} Z(N)\right) = O(\log N).
\end{equation}
Let $1 \le y_{0} < x_{0} \le N$. If $2y_{0} > x_{0}$, then we have
\[
\mathbb{E}\left(\frac{\partial^{2}}{\partial t_{x_{0}} \partial t_{y_{0}}} Z(N)\right) \le \sum_{\substack{(x_{3},x_{4})\\x_{0}+y_{0}=x_{3}+x_{4}\\x_{3} < x_{4}}}\frac{c}{x_{3}^{2/3}}\cdot \frac{c}{x_{4}^{2/3}} + \sum_{\substack{(x_{3},x_{4})\\x_{0}+x_{3}=y_{0}+x_{4}\\1 \le x_{3} < x_{4}\le N}}\frac{c}{x_{3}^{2/3}}\cdot \frac{c}{x_{4}^{2/3}} + \frac{c}{(2x_0 - y_0)^{2/3}}
\]
\[
+ \frac{c}{(2y_0 - x_0)^{2/3}} + \frac{c}{(\frac{x_0 + y_0}{2})^{2/3}}
= \mathbb{E}_{8} + \mathbb{E}_{9} + \mathbb{E}_{10} + \mathbb{E}_{11} + \mathbb{E}_{12}.
\]
In $\mathbb{E}_{8}$, we have $x_{4} > \frac{x_{0}+y_{0}}{2} > \frac{x_{0}}{4}$, and so $\frac{1}{x_{4}^{2/3}} \le \frac{1}{\left(\frac{x_{0}}{4}\right)^{2/3}} = \frac{4^{2/3}}{x_{0}^{2/3}}$, which implies that
\[
\mathbb{E}_{8} = O\left(\sum_{x_{3}=1}^{x_{0}}\frac{1}{x_{3}^{2/3}}\cdot \frac{1}{x_{0}^{2/3}}\right) = O\left(\frac{1}{x_{0}^{1/3}}\right) = O(1). 
\]
Furthermore,
\[
\mathbb{E}_{9} = O\left(\sum_{x_{3}=1}^{\infty}\frac{1}{x_{3}^{2/3}}\cdot \frac{1}{(x_{3}+x_{0}-y_{0})^{2/3}}\right) = O\left(\sum_{x_{3}=1}^{\infty}\frac{1}{x_{3}^{4/3}}\right) = O(1). 
\]
Clearly we have $\mathbb{E}_{10} = O(1)$, $\mathbb{E}_{11} = O(1)$, 
$\mathbb{E}_{12} = O(1)$.

If $2y_{0} \le x_{0}$, then we have
\[
\mathbb{E}\left(\frac{\partial^{2}}{\partial t_{x_{0}} \partial t_{y_{0}}} Z(N)\right) \le \sum_{\substack{(x_{3},x_{4})\\x_{0}+y_{0}=x_{3}+x_{4}\\x_{3} < x_{4}}}\frac{c}{x_{3}^{2/3}}\cdot \frac{c}{x_{4}^{2/3}} + \sum_{\substack{(x_{3},x_{4})\\x_{0}+x_{3}=y_{0}+x_{4}\\1 \le x_{3} < x_{4}\le N}}\frac{c}{x_{3}^{2/3}}\cdot \frac{c}{x_{4}^{2/3}}
\]
\[
+ \frac{c}{(\frac{x_0 + y_0}{2})^{2/3}} + \frac{c}{(2x_0 - y_0)^{2/3}}.
\]
The same computation as in the previous case shows that this is $O(1)$.

\noindent In summary, we have 
\begin{equation}
\mathbb{E}\left(\frac{\partial^{2}}{\partial t_{x_{0}} \partial t_{y_{0}}} Z(N)\right) = O(1).
\end{equation}
Let $1 \le z_{0} < y_{0} < x_{0} \le N$. If $y_{0} + z_{0} > x_{0}$, then
\[
\mathbb{E}\left(\frac{\partial^{3}}{\partial t_{x_{0}} \partial t_{y_{0}} \partial t_{z_{0}} } Z(N)\right) \le \frac{c}{(y_{0}+z_{0}-x_{0})^{2/3}} + \frac{c}{(x_{0}+z_{0}-y_{0})^{2/3}} + \frac{c}{(x_{0}+y_{0}-z_{0})^{2/3}} + 1 
\]
\[
= O(1).
\]
If $y_{0} + z_{0} \le x_{0}$, then
\[
\mathbb{E}\left(\frac{\partial^{3}}{\partial t_{x_{0}} \partial t_{y_{0}} \partial t_{z_{0}}} Z(N)\right) \le \frac{c}{(x_{0}+y_{0}-z_{0})^{2/3}} + \frac{c}{(x_{0}+z_{0}-y_{0})^{2/3}} + 1  = O(1).
\]
In summary, we have
\begin{equation}
\mathbb{E}\left(\frac{\partial^{3}}{\partial t_{x_{0}} \partial t_{y_{0}} \partial t_{z_{0}}} Z(N)\right) = O(1).
\end{equation}
Let $1 \le u_{0} < z_{0} < y_{0} < x_{0} \le N$. If $x_{0} + u_{0} = y_{0} + z_{0}$, then
\[
\mathbb{E}\left(\frac{\partial^{4}}{\partial t_{x_{0}} \partial t_{y_{0}} \partial t_{z_{0}} \partial t_{u_{0}}}Z(N)\right) = 1,
\]
otherwise
\[
\mathbb{E}\left(\frac{\partial^{4}}{\partial t_{x_{0}} \partial t_{y_{0}} \partial t_{z_{0}} \partial t_{u_{0}}}Z(N)\right) = 0.
\]
Then we have
\begin{equation}
\mathbb{E}\left(\frac{\partial^{4}}{\partial t_{x_{0}} \partial t_{y_{0}} \partial t_{z_{0}} \partial t_{u_{0}}}Z(N)\right) = O(1).
\end{equation}
If $\alpha = (\alpha_{1}, \dots{} ,\alpha_{N})$, $\alpha_{i} \ge 0$, $\alpha_{i} 
\in \mathbb{Z}$, $|\alpha| \ge 5$, then clearly we have
\begin{equation}
\mathbb{E}\left(\left(\frac{\partial}{\partial t_{1}}\right)^{\alpha_{1}} \cdots{} \left(\frac{\partial}{\partial t_{N}}\right)^{\alpha_{N}}Z(N)\right) = 0.
\end{equation}
By (5), (10), (11), (12), (13) and (14), we get (8). On the other hand, (10), (11), (12), (13) and (14) proves (9).

By the Kim-Vu inequality with $\lambda = 20\log N$, we get that
\[
\mathbb{P}(|Z(N) - \mathbb{E}(Z(N))| \ge c_{4}(20\log N)^{3.5}\sqrt{\mathbb{E}_{\ge 0}(Z(N))\cdot \mathbb{E}_{\ge 1}(Z(N))}) = O\left(e^{-\frac{20\log N}{4}+3\log N} \right) 
\]
\[
= O\left(\frac{1}{N^{2}}\right).
\] 
Then by the Borel-Cantelli lemma, we get that with probability 1, 
\[
|Z(N) - \mathbb{E}(Z(N))| \le c_{4}(20\log N)^{3.5}\sqrt{\mathbb{E}_{\ge 0}(Z(N))\cdot \mathbb{E}_{\ge 1}(Z(N))},
\]
for every large enough $N$. By (8) and (9), we get (6).

\bigskip
\bigskip

\noindent \textbf{Acknowledgement}: The authors would like to thank Professor
Robert Horv\'ath for the valuable discussions about the numerical estimation of the triple integral.

\end{document}